\documentclass[12pt]{amsart}
\usepackage{amsmath, epsfig, subfigure, afterpage}
\theoremstyle{definition}

\newtheorem{thm}{Theorem}
\newtheorem{cor}[thm]{Corollary}
\newtheorem{lem}{Lemma}
\theoremstyle{remark}

\theoremstyle{plain}

\newcommand{\Z}{{\mathbb Z}}

\newcommand{\R}{{\mathbb R}}
\newcommand{\C}{{\mathbb C}}
\newcommand{\T}{{\mathbb T}}

\newcommand{\N}{{\mathbb N}}

\renewcommand{\a}{\ensuremath{\alpha}}
\newcommand{\g}{\ensuremath{\gamma}}
\renewcommand{\b}{\ensuremath{\beta}}
\newcommand{\del}{\ensuremath{\delta}}

\newcommand{\pfrac}[2]{\left(\frac{#1}{#2}\right)}
\newcommand{\ds}{\displaystyle}
\def\({\left(}
\def\){\right)}

\newcommand{\be}{\begin{equation}}
\newcommand{\ee}{\end{equation}}
\newcommand{\benn}{\begin{equation*}}   
\newcommand{\eenn}{\end{equation*}}

\numberwithin{equation}{section}

\begin{document}

\title[Imaginary parts of zeros of the zeta function]
{On the distribution of imaginary parts of zeros
of the Riemann zeta function}

\author{Kevin Ford and Alexandru Zaharescu}
\address{Department of Mathematics, 1409 West Green Street, University
of Illinois at Urbana-Champaign, Urbana, IL 61801, USA}
\email{KF: ford@math.uiuc.edu, AZ: zaharesc@math.uiuc.edu}


\date{\today}

\thanks{First author supported by National Science Foundation Grants
DMS-0196551 and DMS-0301083.}

\begin{abstract}
We investigate the distribution of the fractional parts of $\a \g$, where
$\a$ is a fixed non-zero real number and $\g$ runs over the imaginary
parts of the non-trivial zeros of the Riemann zeta function.
\end{abstract}

\maketitle

%
%
\section{Introduction}
%

There is an intimate connection between the distribution of the
 nontrivial zeros of the Riemann zeta function $\zeta(s)$ and the
distribution of prime numbers.  Critical to many prime number problems
is the {\it horizontal} distribution of zeros; here the
 Riemann Hypothesis (RH) asserts that the zeros
all have real part $\frac12$.  
There is also much interest in studying the distribution of the
 imaginary parts of the zeros (the {\it vertical} distribution).
For example, one expects
that their consecutive spacings follow the GUE distribution
from random matrix theory.  Originally discovered by
Montgomery \cite{M}, who studied the pair correlation
of zeros of the zeta function,
this phenomenon has been investigated, for higher correlations
and also for more general $L-$functions, by a number of authors,
including Odlyzko \cite{O}, Hejhal \cite{H}, 
Rudnick and Sarnak \cite{RS}, Katz and Sarnak\cite{KS},
Murty and Perelli \cite{MP}, and Murty and Zaharescu \cite{MZ}.

Let $\{ y\}$ denote the fractional part of $y$, which can be
interpreted as the image of $y$ in the torus $\T=\R/\Z$.
In this paper we look at the distribution of $\{ \alpha \gamma\}$
where $\alpha$ is a fixed nonzero real number and $\gamma$
runs over the imaginary parts of the zeros of $\zeta(s)$.
The starting point is
an old formula of Landau \cite{L}, which states that for each fixed $x>1$,
\be\label{eq:Landau}
\sum_{0<\gamma\le T} x^{\rho} = -\frac{T}{2\pi} \Lambda(x) + O(\log T).
\ee
Here the sum is over zeros $\rho=\b + i\gamma$ of $\zeta$, 
$\Lambda(x)$ is the von Mangoldt function for integral $x>1$ and
$\Lambda(x)=0$ for non-integral $x>1$.  
Put $x=e^{2\pi j \alpha}$ into \eqref{eq:Landau}, where $\alpha>0$ and
$j$ is a positive integer.
If RH holds, then \eqref{eq:Landau} implies that
\be
\sum_{0<\gamma\le T} e^{2\pi i j \a \gamma} = O(T),
\ee
while there are $\gg T\log T$ summands on the left side (see \eqref{NT}
below).  Thus, as Rademacher observed in 1956 \cite{R},
RH plus Weyl's criterion (\cite{KR}, Ch. 1, Theorem 2.1) 
implies that the sequence
$\{\alpha \gamma\}$ is uniformly distributed in $[0,1)$.
In 1975, Hlawka \cite{Hl} proved this conclusion unconditionally.
His result depended on a version of \eqref{eq:Landau} which is uniform in
$x$ and results about the density of zeros lying off the line $\beta=\frac12$
(see Lemmas 1,2 below).
Hlawka also showed via the Erd\H os-Tur\'an inequality (\cite{KR}, Ch. 2,
Theorem 2.5) that the discrepancy of the set
$\{ \{\a\g\}: 0<\g\leq T\}$ is $O(\frac{1}{\log T})$ assuming RH.
More recently Fujii \cite{F1} showed unconditionally that this discrepancy
is $O(\frac{\log\log T}{\log T})$.

Ramemacher also noted in \cite{R} that if $\a = \frac{k\log p}{2\pi}$, where
$p$ is a prime and $k\in \N$ (i.e. $x=p^k$ so $\Lambda(x)=\log p$),
by \eqref{eq:Landau} there should be a ``predominance of
terms which fulfill $|\{\a\g\} - 1/2| < 1/4$''.  
We will give a very precise meaning to this statement, and also
show that Hlawka's discrepancy result is best possible in a certain sense.
We will see that there is a certain limiting measure
$\mu_{\alpha}$
that one can naturally associate to each $\alpha>0$,
which explains among other things the behavior
of the above discrepancy as $T\rightarrow\infty.$

For any real numbers $\alpha, T >0$ 
consider the measure $\mu_{\alpha, T}$ defined on $\T$ by
\begin{equation}\label{b1}
\mu_{\alpha ,T} = \frac 1T \sum_{0 < \gamma \le T} 
\delta_{\{\alpha\gamma\}} - \frac {N(T)}T \mu.
\end{equation}
Here $\mu$ denotes the Haar measure on $\T,$
$\delta_{\{\alpha \gamma\}}$ is a unit point
delta mass at $\{\alpha \gamma\}$ and $N(T)$ is the number of zeros 
$\rho=\beta+i\gamma$
of $\zeta$, counted with multiplicity, with $0<\beta<1$ and $0<\gamma
\le T$.
To the measure $\mu_{\a,T}$ we associate the function
\be\label{Mdef}
M(y;T)=M_\a(y;T) := \mu_{\a,T}([0,y))=
\frac{1}{T} \sum_{\substack{0<\g \le T \\ \{\a\g\} < y}}
1 - y \frac{N(T)}{T}. 
\ee
In particular, $M(0;T)=M(1;T)=0$.  Let
\be\label{DM}
D_\a^*(T) = D^*(T)= \sup_{0\le y \le 1} \frac{T}{N(T)} |M(y;T)|
\ee
denote the discrepancy of the set $\{ \{\a\g\}:0<\g\le T\}$.

For fixed $\alpha$, the above measures $\mu_{\alpha, T}$, seen
as continuous linear functionals on $C(\T),$ form an unbounded set.
In fact, the total variation of $\mu_{\alpha,T}$
grows like $\log T$ as $T\rightarrow\infty$, since
(\cite{T}, Theorem 9.4)
\be\label{NT}
N(T) = \frac{T}{2\pi} \log \frac{T}{2\pi} - \frac{T}{2\pi} + O(\log T)
\quad (T\ge 1).
\ee
We will see, however, that for a large class of
functions $h:\T\rightarrow\C$ (including all absolutely continuous
functions on RH),
$\lim_{T\rightarrow\infty}\int_{\T}hd\mu_{\alpha,T}$ 
exists.  Moreover, there is a  unique Borel measure 
$\mu_{\alpha}$ on $\T$
which is absolutely continuous with respect to the
Haar measure, such that for such $h$ one has
\begin{equation}\label{b2}
\lim_{T\rightarrow\infty}\int_{\T}hd\mu_{\alpha,T}=
\int_{\T}hd\mu_{\alpha}.
\end{equation}
Let $g_{\alpha}:[0,1)\rightarrow\R$ denote,
via the above identification between $\T$ and $[0,1),$ 
the density of this measure.
Note that $\mu_{\alpha,T}(\T)=0$ for any $T,$ hence

\begin{equation}\label{b3}
\int_0^1g_{\alpha}(t)dt=0.
\end{equation}

We have identified the density function
$g_{\alpha},$ and we found that
\be\label{b4}
g_{\alpha}(t)=0 \quad (t\in\T)
\ee
provided $\alpha$ is not a rational multiple of a
number of the form $\frac {\log p}{\pi}$ with $p$ prime. 
If $\alpha= \frac {a\log p}{2\pi q}$ for some prime number $p$
and positive integers $a,q$ with
$(a,q)=1,$ then for any $t \in [0,1),$ 
\begin{equation}\label{b5}
\begin{split}
g_{\alpha}(t) &= -\frac{\log p}{\pi} \Re \sum_{k=1}^\infty (p^{a/2} 
e^{2\pi i qt})^{-k} \\
&=-\frac{(p^{\frac a2}\cos 2\pi qt -1) \log p } 
{\pi(p^a-2p^{\frac a2}\cos 2\pi qt+1)}\,. 
\end{split}
\end{equation}

As a function of $t$, $g_{\alpha}(t)$ attains its
global minimum at each of the points $t=\frac kq$, $k = 0,1,\cdots, q-1,$
the minimum being 
\begin{equation}\label{b8}
g_{\alpha}\Big (\frac kq\Big )=-\frac {\log
p}{\pi (p^{\frac a2} -1)} <0. 
\end{equation}
In particular, this shows that there is a shortage of
zeros of  $\zeta(s)$ with imaginary parts $\gamma$ such that $\{\alpha
\gamma\} =\{\frac {a\log p}{2\pi q}\gamma\}$ is close to one of the
points $\frac kq$, $k = 0,1,\ldots, q-1$.  When $q=1$, this 
corresponds to Rademacher's statement mentioned above.

We conclude this introduction by showing some histograms of $M(y;T)$
for $T=600,000$ ($N(T)=999508$)
 and a few values of $\a$.  The list of zeros to
this height were kindly supplied by Andrew Odlyzko.
We partition $[0,1)$ into 500 subintervals of length $\frac{1}{500}$.
In Figure 1 we plot for each subinterval $I=[y,y+\frac{1}{500})$ the value
of $500 (M(y+\frac{1}{500})-M(y))$ and also the graph of $g_\a(y)$.

%
%

\afterpage{\clearpage}
\begin{figure}[t]
\parbox{.8\linewidth}{\subfigure{
\epsfig{file=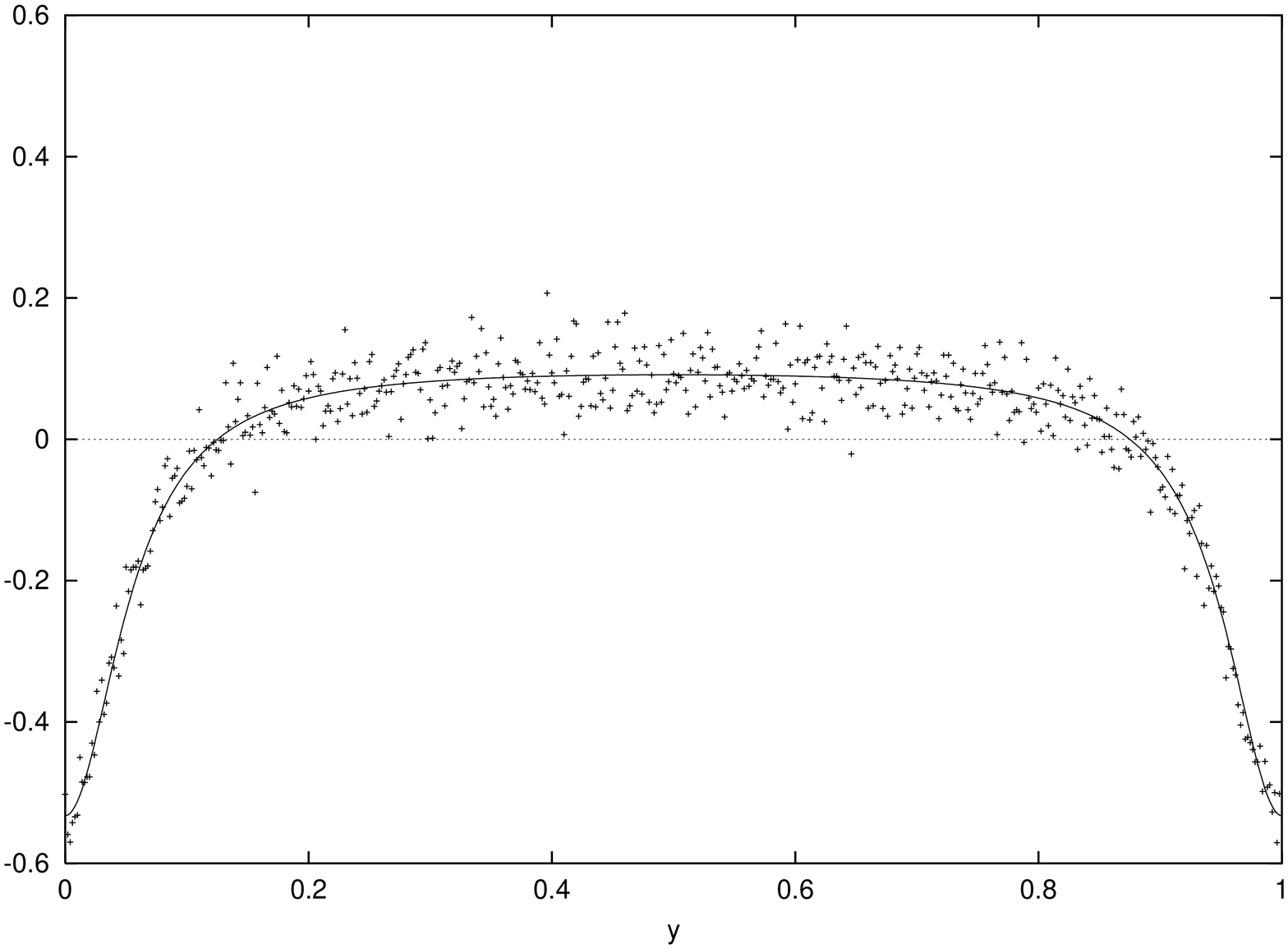, height=4in, width=2.2in, angle=270}}}
\hfill  $\ds \a=\frac{\log 2}{2\pi}$ \hfill
\parbox{.8\linewidth}{\subfigure{
\epsfig{file=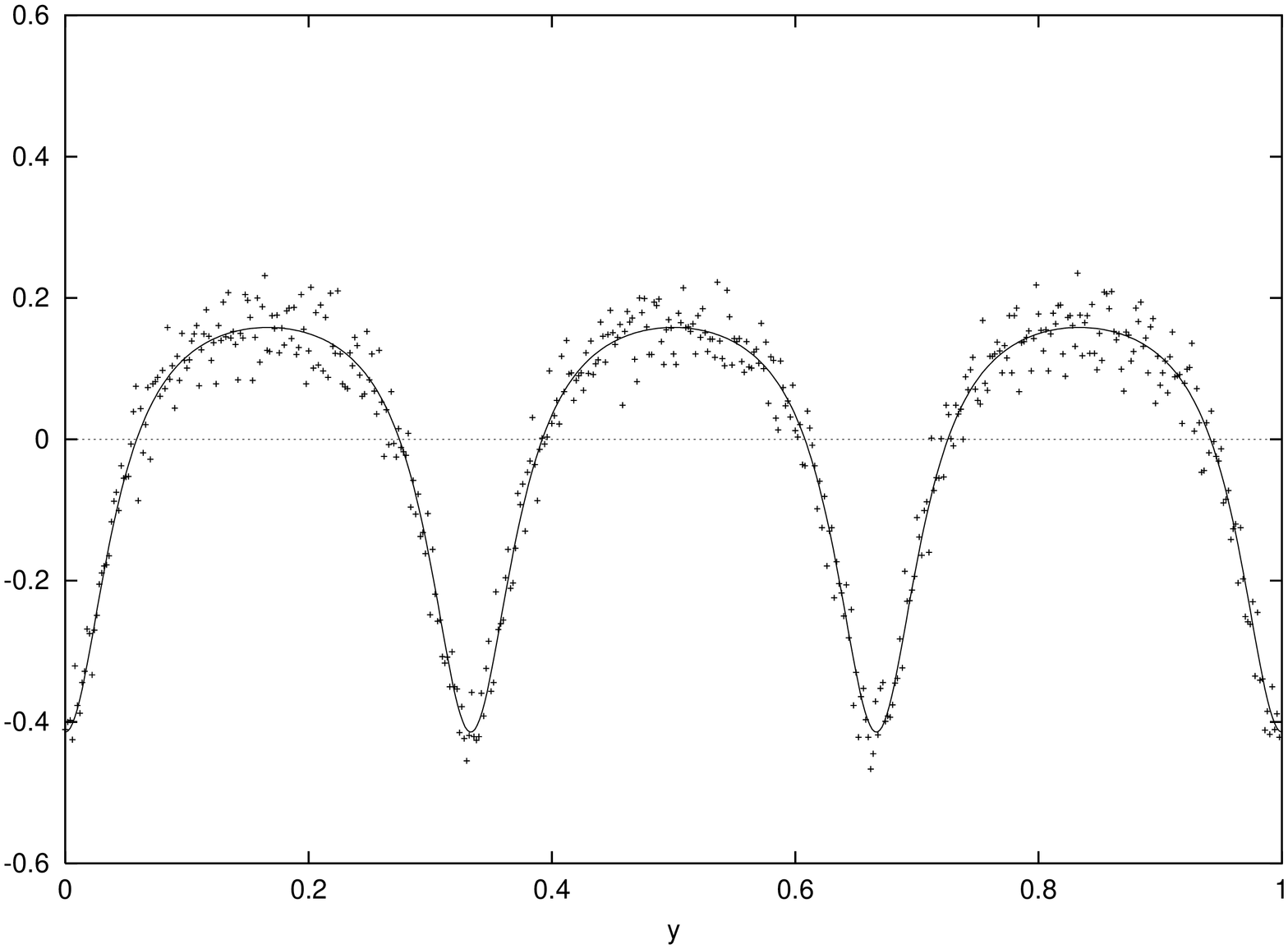, height=4in, width=2.2in, angle=270}}}
\hfill  $\ds \a=\frac{\log 5}{3\cdot 2\pi}$ \hfill
\parbox{.8\linewidth}{\subfigure{
\epsfig{file=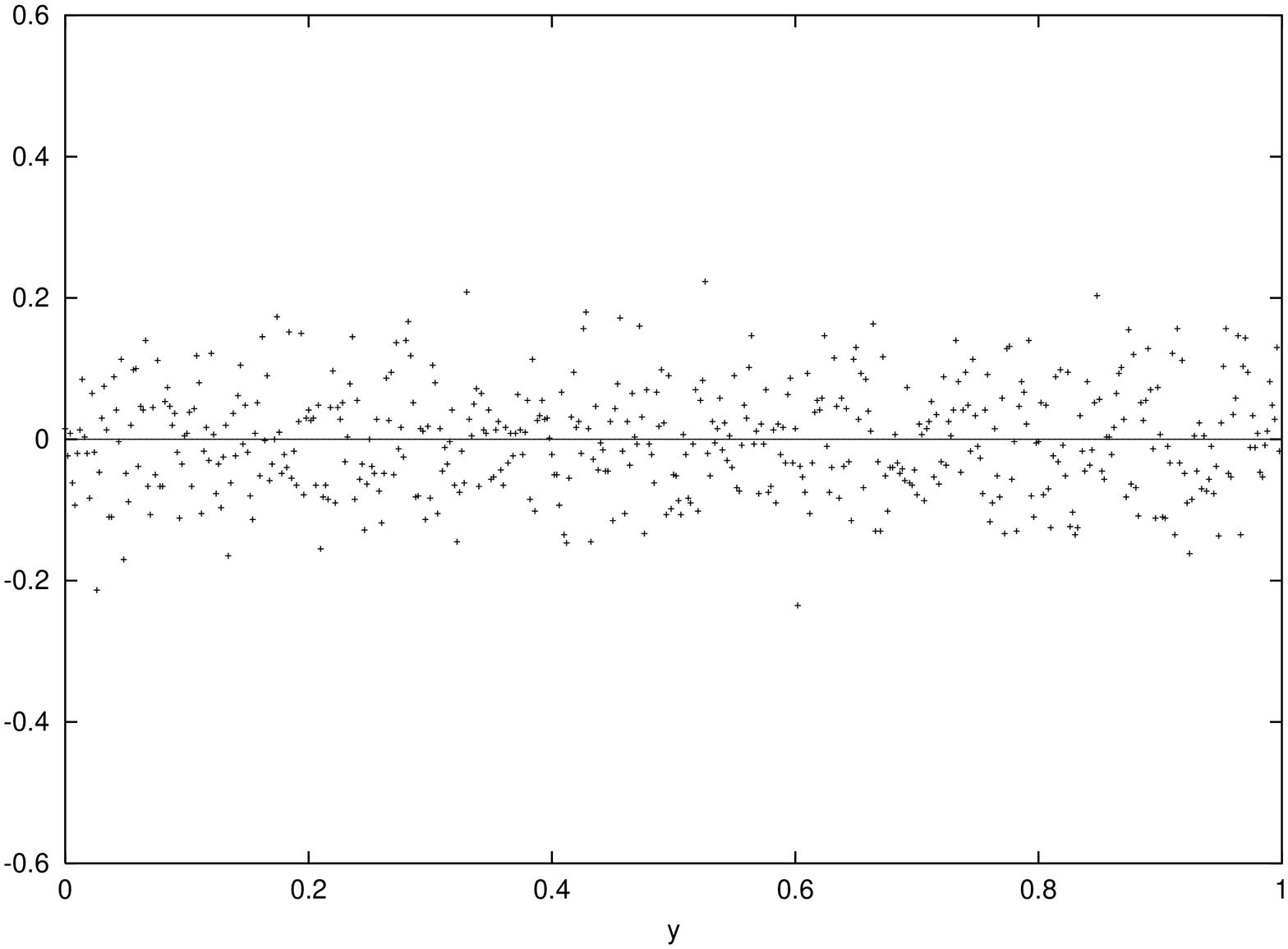, height=4in, width=2.2in, angle=270}}}
\hfill  $\ds \a=\frac{\log 6}{2\pi}$ \hfill
\caption{$500(M(y+\frac{1}{500})-M(y))$ vs. $g_\a(y)$ for
$T=600000$.}
\end{figure}

%
\section{Statement of results and conjectures}
%

To state our main theorem, we consider for any $f\in L_1(\T)$ the
modulus of continuity (\cite{BN}, Definition 1.5.1)
\benn
\omega(f;\delta) = \sup_{|u| \le \delta} \int_\T |f(t+u)-f(t)|\, dt.
\eenn
For any real numbers $\alpha>0$ and $T>0$ define the
measure $\mu_{\alpha,T}$ as in \eqref{b1},
and let $g_{\alpha}:\T\rightarrow\R$ be defined by \eqref{b4}
or \eqref{b5}, as appropriate. 

%
%

\begin{thm}\label{thm1}
Let $\a > 0$ and let
$h:\T\rightarrow\C$ be an absolutely continuous function with
\be\label{thm1cond}
\omega(h';\delta) = o(1/\log (1/\del)) \quad (\del\to 0^+).
\ee
Then
\be\label{conv}
\lim_{T\rightarrow\infty}\int_{\T}hd
\mu_{\alpha,T}=\int_{\T}h
g_{\alpha}d\mu.
\ee
\end{thm}

%
%

\begin{cor}\label{cor2}
The equality \eqref{conv} holds for all $h\in C^2(\T)$.
\end{cor}

\begin{proof}
Since $h\in C^2(\T)$,
\benn
\omega(h';\del) = \sup_{|u| \le \del} \int_\T \left| 
\int_{t}^{t+u} h''(y)\, dy\right|\, dt \le \del \int_\T |h''| = O(\del).
\eenn
\end{proof}

\begin{cor}\label{cor4}
Suppose $\a=\frac{a\log p}{2\pi q}$, where $p$ is prime and $(a,q)=1$.
Then
\benn
\int_\T |M(y;T)| \, dt \ge B(\a) + o(1),
\eenn
where
\benn
B(\a) = \frac{\int_\T g_\a^2\, d\mu}{\max |g'_\a|}.
\eenn
Therefore, $D_\a^*(T) \ge \frac{B(\a)/2\pi+o(1)}{\log T}$ for these $\a$.
\end{cor}

\begin{proof}
Apply Theorem \ref{thm1} with $h=g_\a$ (this may or may not be optimal).
By integration by parts,
\benn
\begin{split}
\left| \int_\T g_\a^2 + o(1) \right| &= \left| \int_\T M(y;T) g_\a'(y)\, dy
\right| \\
&\le \(\int_\T |M(y;T)|\, dy \) \max |g_\a'(y)|.
\end{split}
\eenn
\end{proof}

%
%

\begin{thm}\label{thm2}
Let $\a>0$ and suppose $D_\a^*(T) \ll \frac{1}{\log T}$.  
Then \eqref{conv} holds for all
absolutely continuous functions $h:\T\rightarrow\C$.
\end{thm}

\begin{cor}\label{cor3}
On RH, \eqref{conv} holds for all absolutely continuous functions 
$h$ on $\T$.
\end{cor}

\begin{proof}
By Hlawka's Theorem, if RH is true then $D_\a^*(T) \ll 1/\log T$.
\end{proof}

%
%

The main open problem in this line of investigation is to determine the
largest ``natural'' class of functions for which \eqref{conv} holds.
We conjecture that \eqref{conv} holds for the characteristic function
of the interval $[0,y)$.

\noindent
{\bf Conjecture A.}
If $\a=\frac{a\log p}{2\pi q}$, where $p$ is prime and $(a,q)=1$, then
uniformly in $y$,
\benn
\lim_{T\to \infty} 
M(y;T) = \int_{0}^y g_\a(t)\, dt = -\frac{\log p}{2\pi^2 q}
\arg (1-p^{-a/2} e^{-2\pi i q y}).
\eenn
When $\a$ is not of this form, $\lim_{T\to\infty} M(y,T) = 0$ uniformly
in $y$.
\medskip


\begin{cor}\label{conj3}  
Assume Conjecture A.
Suppose $\a=\frac{a\log p}{2\pi q}$, where $p$ is prime and $(a,q)=1$.
As $T\to\infty$,
\benn
D_\a^*(T) =(1+o(1)) \frac{\log p}{\pi q}\; \frac{\arcsin(p^{-a/2})}{\log T}.
\eenn
When $\a$ is not of this form, $D_\a^*(T)=o(1/\log T)$.
\end{cor}

\bigskip

In the opposite direction, \eqref{conv} cannot hold for all functions
$h$ which are continuous and differentiable on $\T$.  This is a 
consequence of a property of general sequences (uniformly distributed
or not), which we state below.

\begin{thm}\label{negative}
Let $a_1, a_2, \ldots$ be an arbitrary sequence of numbers in $\T$, let
$t$ be a point in $\T$ and let $f(x)$ be a function decreasing monotonically
to $0$ arbitrarily slowly.  Then there is a function $h$, continuous
and differentiable on $\T$, and which is $C^\infty(\T \backslash \{ t \})$,
so that for an infinite set of $n\in \N$,
\begin{equation}\label{eq2.3}
\left| \frac{1}{n} \sum_{j=1}^n h(a_j) - \int_\T h \right| \ge f(n).
\end{equation}
\end{thm}

In connection with his work
on $k$-functions,  Kaczorowski \cite{K} was also led to study
similar questions on the distribution of the imaginary parts
of zeros of the zeta function and Dirichlet $L$-functions.  
In our terminology, his results concern the measures
\benn
\nu_{\a,n} = \frac{1}{n!} \sum_{\g>0}  e^{-\g}
\g^n \delta_{\{\a \g\}} - \frac{1}{n!} \sum_{\g>0} e^{-\g} \g^n \mu,
\eenn
which essentially capture the distribution of $\{ \a\g\}$ for
$|\g-n| \ll \sqrt{n}$.
Kaczorowski's Conjecture A corresponds to our
Corollary \ref{conj3}, and his Conjectures B, B$_1$ and B$_2$ to
our Conjecture A.

We also mention that Fujii (\cite{F2}, Theorem 5) 
has proved \eqref{conv} in the
special case $h(u)=B_2(u)=u^2-u+\frac16$ (which is covered by Theorem 1)
He has also proven, for a wide class of smooth functions $f$, 
that the sequence $\{ f(\gamma) \}$ is uniformly distributed \cite{F3}.

Lastly, we mention that analogs of our results should hold
for a wide variety of $L$-functions, including Dirichlet $L$-functions.
In deriving such results,
the primary tools would be a generalization of Lemma 1 below
and zero density estimates
analogous to Lemma 2 below which are nontrivial for certain sequences of
$\sigma, T$ with $T\to\infty$ and $\sigma - \frac12 \to 0^+$.
This will be addressed in a future paper.
\bigskip

%
\section{Proof of Theorem \ref{thm1}}
%

Our principal tools are the following two lemmas.

\begin{lem}\label{lem1} Let $x,T>1,$ and denote by $n_x$ the 
nearest prime power to $x.$ Then
\be\label{c1}
\sum_{0<\gamma\leq T}x^{\rho}=
-\frac{\Lambda(n_x)}{2\pi} \frac{e^{iT\log(x/n_x)}-1}
{i\log(x/n_x)}
+O\(x\log^2(2xT) + \frac{\log 2T}{\log x}\),
\ee
where if $x=n_x$ the first term is $-T \frac{\Lambda(n_x)}{2\pi}$.
\end{lem}
This is a uniform version of a theorem of Landau \cite{L}, and
the proof is nearly identical to the proof of Theorem 1 of Gonek (\cite{G},
\S 3).  The only difference is in the treatment of the 
term $n=n_x$ occuring in Gonek's integral $I_1$.  That term is 
\benn
V=-\frac{\Lambda(n_x)}{2\pi} (x/n_x)^c \int_1^T (x/n_x)^{it}\, dt,
\quad c=1+\frac{1}{\log(3x)}.
\eenn
Gonek estimates this crudely as $V \ll \log x \min (T, x/|x-n_x|)$,
but we are more precise.
If $|x-n_x| \ge 1$, then $V\ll x\log(2x)$ and otherwise
\benn
\begin{split}
V &= -\frac{\Lambda(n_x)}{2\pi} (1 + O(|x-n_x|/x)) 
\( \frac{e^{iT\log(x/n_x)}-1}
{i\log(x/n_x)} + O(1)\) \\
&= O(\log(2x)) -\frac{\Lambda(n_x)}{2\pi} \frac{e^{iT\log(x/n_x)}-1}
{i\log(x/n_x)}.
\end{split}
\eenn

Let $N(\sigma,T)$ be the number of zeros $\rho=\beta+i\g$ of $\zeta(s)$
with $|\g| \le T$ and $\b \ge \sigma$.  The following density bound
 is due to Selberg (Theorem 9.19C of \cite{T}).

\begin{lem}\label{lem2}
Uniformly in $\frac12\leq\sigma\leq1$, we have
$N(\sigma,T)\ll T^{1-\frac14 (\sigma-1/2)}\log T$.
\end{lem}

We now proceed with the proof of Theorems \ref{thm1} and \ref{thm2}.
Fix $\alpha>0$ and an absolutely continuous function
$h:\T\rightarrow\C$.  In what follows, constants implied by the
$O-$ and $\ll-$symbols may depend on $\a$ and $h$. 
Suppose either $D_\a^*(T)\ll 1/\log T$ or $h$ satisfies \eqref{thm1cond}.
Looking at the definition
of $\mu_{\alpha,T}$ we see that what we need to 
show in order to complete the proof of the theorems
is that
\be\label{main}
\frac1T
\sum_{0 < \gamma\leq T}h(\alpha\gamma)-\frac{N(T)}T
\int_0^1 h(u)du=\int_0^1 h(u)g_{\alpha}(u)\, du+o(1)
\ee
as $T\rightarrow\infty.$
By treating separately the real part and
the imaginary part of $h,$ we may assume in 
what follows that $h$ is real.
We start by approximating $h$ by a trigonometric polynomial based on
its Fourier series
\be\label{c4}
h(u)=\sum_{m\in\Z}c_m e^{2\pi mui},
\ee
where the Fourier coefficients $c_m$ are given by
\be\label{c5}
c_m=\int_0^1 h(u) e^{-2\pi mui}\, du.
\ee
Since $h$ is absolutely continuous, we have
$h'\in L_1(\T)$
and consequently by the Riemann-Lebesgue Lemma,
\be\label{cm}
c_m = o(1/m) \qquad (m\to \infty).
\ee

Let
$$
K_J(u) = \frac{3}{J(2J^2+1)} \pfrac{\sin \pi J u}{\sin \pi u}^4 = 
\frac{3}{J(2J^2+1)} \biggl( \sum_{k=0}^{J-1} e^{(1-J+2k)u\pi i}
\biggr)^4
$$
be the Jackson  kernel (\cite{BN}, Problem 1.3.9) and define
\benn
H_J(y) = \int_{\T} K_J(u) h(u-y)\, du = 
\sum_{|j| \le 2J} A^{(J)}_j c_j e^{2\pi ju i},
\eenn
where
$$
A^{(J)}_0=1, \quad |A^{(J)}_j| \le 1, \quad \lim_{J\to\infty}
A^{(J)}_j=1\text{ for each fixed }j.
$$
This kernel is chosen because it provides a very fast rate 
of $L_1(\T)$-convergence of $H_J$ to $h$ (see \cite{BN}, \S 1.6 and
Ch. 2 for more on this subject).
The same is true for the convergence of $H_J'$ to $h'$. 
Specifically, we have (\cite{BN}, Lemma 1.5.4, Theorem 1.5.8 and
Corollary 1.5.9)
\be\label{rate}
\int_\T |H_J'-h'| = O(\omega(h';1/J)).
\ee

Let $T$ be a large real number.
Applying integration by parts to the left side of \eqref{conv} and using 
\eqref{Mdef} and \eqref{DM}, we obtain
\be
\begin{split}\label{main1}
\int_{\T}hd\mu_{\alpha,T} &= \int_\T h(y)\, dM(y;T) \\
&=\int_\T H_J(y) \, dM(y;T) + \int_\T (h(y)-H_J(y))  \, dM(y;T)\\
&=\int_\T H_J(y) \, dM(y;T) - \int_\T M(y;T)\, (h'(y)-H_J'(y))\, dy \\
&=\frac{1}{T} \sum_{1\le |j| \le 2J-2} A^{(J)}_j  c_j 
\sum_{0<\g\le T} x_j^{i\g} + O\(D^*(T) \log T \int_\T |h'-H_J'|\),
\end{split}
\ee
where we have written $x_j=e^{2\pi j \a}$.

To apply Lemma \ref{lem1}, we write
$x^{i\g} = x^{\rho-1/2} + x^{i\g} (1-x^{\beta-1/2})$.
Let $\delta=50 \frac{\log\log T}{\log T}$.
Assume that 
$$
0<\delta\log x < 1.
$$
Note that if $\beta+i\g$ is a zero, then $1-\beta+i\g$ is also a zero.
We obtain from Lemma \ref{lem2} the bounds
\begin{align*}
\sum_{\substack{0<\g \le T \\ |\beta - \frac12| \ge \delta}}  
  &x^{i\g} (1-x^{\beta-1/2}) \ll \sum_{\substack{0<\g \le T \\ \beta \ge 
  \frac12+\delta}} x^{\beta-1/2} \\
&\ll x^{\delta} N(1/2+\delta,T)+ \log x \int_{1/2+\delta}^1 
  x^{\sigma-1/2} N(\sigma,T)\, d\sigma \\
&\ll T \log T \biggl[ T^{-\del/4} + \log x  \int_{\del}^{1/2}
  \pfrac{T^{1/4}}{x}^\theta \, d\theta \biggr] \\
&\ll \frac{T}{\log^{10} T}
\end{align*}
and
\begin{align*}
\sum_{\substack{0<\g \le T \\ |\beta - \frac12| < \delta}}  
   x^{i\g} (1-x^{\beta-1/2}) &\ll \sum_{\substack{0<\g \le T \\ \frac12 < \beta
   < \frac12+\delta}} (x^{\beta-1/2}+x^{1/2-\beta}-2) \\
&\ll  \log^2 x \sum_{\substack{0<\g \le T \\ \frac12 < \beta
   < \frac12+\delta}} (\beta-1/2)^2 \\
&\ll \log^2 x \int_{0}^\delta \theta N(1/2+\theta,T)\, d\theta \\
&\ll \frac{T\log^2 x}{\log T}.
\end{align*}
Thus
\be\label{sumzeros}
\sum_{0<\g \le T} (x^{i\g} - x^{\rho-1/2}) \ll \frac{T\log^2 x}{\log T}.
\ee
We apply these estimates with $x=x_j$, $1\le j\le 2J-2$, noting that 
$x_{-j}^{i\g} = x_j^{-i\g}$ and $c_{-j}=\overline{c_j}$.
By Lemma \ref{lem1}, we obtain
\be\label{main2}
\begin{split}
  \int_{\T}&hd\mu_{\alpha,T} = -\frac{1}{\pi} \Re \left[\sum_{j=1}^{2J-2} 
  \frac{A^{(J)}_j c_j \Lambda(n_{x_j})}{\sqrt{x_j}} 
  \frac{\sin\(T\log (x_j/n_{x_j})\)} {T\log (x_j/n_{x_j})} \right] \\
& \qquad + O \( \frac{1}{T}\sum_{j=1}^{2J-2} |c_j| x_j^{1/2} \log^2 (2x_j T)
  + \frac{1}{\log T} \sum_{j=1}^{2J-2} |c_j| \log^2 x_j\) \\
& \qquad + O\( D^*(T) \log T \int_\T |h'-H_J'| \).
\end{split}
\ee
Now put
\be\label{J}
J = \lfloor (\log T)^{1/3} +1\rfloor.
\ee
Using \eqref{cm}, \eqref{J} and $\log x_j =O(j)$,
we find that the first error term in \eqref{main2} is
\benn
O\( \frac{x_{2J-2}^{1/2} \log^2 T}{T} + 
\frac{J^2}{\log T}\) = o(1) \quad (T\to \infty). 
\eenn
We now estimate the second error term in \eqref{main2}.
Since $h'\in L_1(\T)$, $\int_\T |h'-H_J'| \to 0$ as 
$J\to \infty$.  Thus, if $D_\a^*(T) \ll 1/\log T$, then 
\be\label{hH}
D^*(T) \log T \int_\T |h'-H_J'| = o(1).
\ee
Similarly, \eqref{hH} follows if \eqref{thm1cond} holds by
\eqref{rate}.

Consider now the main term in \eqref{main2}.  This sum is absolutely
and uniformly convergent in $T$ and each term with $n_{x_j} \ne x_j$
tends to zero as $T\to\infty$.  Thus, the aggregate of such terms is
$o(1)$ as $T\to\infty$.  
Therefore, by \eqref{sumzeros} and \eqref{hH},
$$
\int_{\T}hd\mu_{\alpha,T} = -\frac{1}{\pi} \Re \left[\sum_{j=1}^{2J-2} 
  \frac{A^{(J)}_j c_j \Lambda(x_j)}{\sqrt{x_j}} \right] 
+ o(1) \quad (T\to\infty).
$$

In the case where $\alpha$ is not of the form
$\frac aq\frac{\log p}{2\pi}$
for positive integers $a, q$ and prime $p$, $\Lambda(x_j)=0$ for
all $j$ and this finishes the proof.
If $\a$ does have this form, $\Lambda(x_j)=0$ unless $q|j$, in which case
$\Lambda(x_j)=\log p$.  Thus
$$
\int_{\T}hd\mu_{\alpha,T} = -\frac{\log p}{\pi} \Re \left[\sum_{m \le
 (2J-2)/q} 
  \frac{A^{(J)}_{qm} c_{qm}}{p^{am/2}} \right] + o(1) \quad (T\to\infty).
$$
It is easily seen that (since $h$ and $g_\a$ are absolutely continuous)
that the right side approaches $\int h g_\a$ as $T\to\infty$.
This completes the proof of Theorem 1.

%
\section{Proof of Theorem \ref{negative}}
%

First, if $a_1, a_2, \ldots$ is not uniformly distributed in $[0,1)$,
by Weyl's criterion the conclusion of Theorem 6 holds with 
$h(y)=e^{2\pi i k y}$ for some integer $k$, regardless of the
function $f$.

From now on assume that $a_1, a_2, \ldots$ is uniformly distributed 
in $[0,1)$. Let $r(x)$ be a $C^\infty(\R)$ function such that
$r(x)=0$ for $|x| \ge 1/2$, $r(0)=1$ and $r$ is monotone on $[-1/2,0]$
and $[0,1/2]$. For any real number $u$, and any $v>0$, $\delta>0$,
denote by $r_{u,v,\delta}$ the function defined by
\begin{equation*}
r_{u,v,\delta}(x)=vr\left(\frac{x-u}{\delta}\right), \,\,\, x\in\R .
\end{equation*}
Thus $r_{u,v,\delta}$ is nonnegative, is $C^\infty(\R)$, it is
supported on an interval of length $\delta$ centered at $u$, and
one has $r_{u,v,\delta}(u)=v$. Note also that 
\begin{equation*}
\int_{\R}r_{u,v,\delta}(x) dx=v\delta\int_{\R}r(x) dx.
\end{equation*}
We now proceed to construct a function $h$ as in the statement of the
theorem. It is enough to define $h$ on the interval $(t-1,t]$.
We set $h(t)=0$. Next, we write the interval $(t-1,t)$ as
a disjoint union of intervals 
\begin{equation*}
I_k=\left( t-\frac1{2^{k-1}},t-\frac1{2^k}\right],\,\,\, k=1, 2, 3,\dots,
\end{equation*}
and we define $h$ inductively on each of the intervals $I_k$.
We also set
\begin{equation*}
J_k=\left[ t-\frac7{2^{k+2}},t-\frac5{2^{k+2}}\right]
\end{equation*}
for any $k\geq1$. Thus $J_k$ is contained in $I_k$ and 
length $(J_k)=1/2^{k+1}$.
For each $k\geq1$, we will select a positive
integer $n_k$, and two real numbers $v_k>0$, $\delta_k>0$.
Then $h$ will be defined on $I_k$ by
\begin{equation}\label{eq4.1}
h(x)=\sum_{\substack{a_n\in J_k\\n\leq n_k}}r_{a_n,v_k,\delta_k}(x).
\end{equation}
Here the sum is restricted to distinct values of $a_n$. Thus if
$a_n$ is the same for several values of $n\leq n_k$, only one
of these values of $n$ is taken on the right side of \eqref{eq4.1}.
We will choose $\delta_k$ to be smaller than $1/2^{k+2}$. Then
the right side of \eqref{eq4.1} will indeed be supported inside
the interval $I_k$. Moreover, $h$ will vanish on a small interval
around each of the points $t-\frac1{2^k}$ with $k\geq1$.
Note also that on each $I_k$, $h$ is a finite sum of $C^{\infty}$
functions, so after this construction is complete $h$
will be $C^{\infty}$ on $\T\setminus \{t\}$.
For each $k$, we only choose $\delta_k$ after $n_k$ has already
been chosen. Then we let $\delta_k$ to be small enough so
that $|a_n-a_{n'}|>\delta_k$ for any $n,n'\leq n_k$ 
with $a_{n'}\neq a_n$. This will make the supports of
the functions $r_{a_n,v_k,\delta_k}$ from the right side of 
\eqref{eq4.1} to be disjoint. As a consequence, one will have
$0\leq h(x)\leq v_k$ for any $x\in I_k$. If we choose the sequence
$(v_k)_{k\geq1}$ to be decreasing to $0$, the function $h$ will
be continuous at $t$. If furthermore the sequence 
$(v_k)_{k\geq1}$ is chosen so that 
\begin{equation*}
\lim_{k\rightarrow\infty}2^kv_k=0,
\end{equation*}
then $h$ will be differentiable at $t$, and $h'(t)=0$. We put
$v_k=1/3^k$, so the above condition holds.
It remains to construct the sequences $(n_k)_{k\geq1}$,
and $(\delta_k)_{k\geq1}$.
Fix a $k$ and assume that $n_j$ and $\delta_j$ have already
been defined for $j=1,\dots,k-1$. This means that $h$ has 
been constructed on the interval $E_k:=I_1\cup\dots\cup I_{k-1}$.
Let $n_k'$ be such that for any $n\geq n_k'$ we have
$f(n)<1/7^k$.
Now, since the sequence $(a_n)_{n\geq1}$ is uniformly distributed,
there exists an $n_k''$ such that for any $n\geq n_k''$
one has
\begin{equation*}
\left| \frac1n\sum_{\substack{1\leq j\leq n\\
a_j\in E_k}}h(a_j)-\int_{E_k}h\right|\leq \frac1{7^k}.
\end{equation*}
The fact that $(a_n)_{n\geq1}$ is uniformly distributed
also implies the existence of an $n_k'''$ such that for
any $n\geq n_k'''$,
\begin{equation*}
\#\{ m\leq n: a_m\in J_k\} > \frac{n}{2^{k+2}}\, .
\end{equation*}
We put $n_k=\max\{ n_k',n_k'',n_k'''\}$. Then by the above it
follows that
\begin{equation*}
f(n_k)<\frac1{7^k}\, ,
\end{equation*}
\begin{equation*}
\left| \frac1{n_k}\sum_{\substack{1\leq j\leq n_k\\
a_j\in E_k}}h(a_j)-\int_{E_k}h\right|\leq \frac1{7^k}\, ,
\end{equation*}
and
\begin{equation*}
\frac1{n_k}\sum_{\substack{1\leq j\leq n_k\\
a_j\in I_k}}h(a_j)\geq \frac{\#\{ j\leq n_k: a_j\in J_k\}}{n_k}v_k
\geq \frac1{2^{k+2}3^k}\, .
\end{equation*}
Lastly, we choose $\delta_k>0$ to be small enough so that
it satisfies all the requirements stated so far, and
such that we also have
\begin{equation*}
\int_{I_k}h < \frac1{7^k}.
\end{equation*}
This completes the construction of the three sequences
$(n_k)_{k\geq1}$, $(v_k)_{k\geq1}$, $(\delta_k)_{k\geq1}$,
and thus also the construction of $h$.
It remains to check that the inequality from the statement
of the theorem holds for infinitely many $n$. Take $n=n_k$,
and break the interval $(t-1,t]$ as a disjoint union of three
intervals, $(t-1,t]=E_k\cup I_k\cup T_k$, where
$T_k=(t-1/2^k,t]$.
We also break accordingly the sum and the integral 
from \eqref{eq2.3}, 
\begin{equation*}
\frac{1}{n_k} \sum_{j=1}^{n_k} h(a_j)=\Sigma_{E_k}+\Sigma_{I_k}+\Sigma_{T_k},
\end{equation*}
and respectively
\begin{equation*}
\int_\T h =\int_{E_k}h +\int_{I_k}h+\int_{T_k}h. 
\end{equation*}
We know that $\left|\Sigma_{E_k}-\int_{E_k}h\right|\leq 1/7^k$,
$\Sigma_{I_k}\geq 1/(4\cdot 6^k)$,
 $0<\int_{I_k}h\leq 1/7^k$ and $f(n_k)\leq1/7^k$.
Also $\Sigma_{T_k}\geq0$, and since for any $m>k$ we have 
$0\leq\int_{I_m}h\leq1/7^m$, it follows that $\int_{T_k}h<2/7^k$.
We conclude that for large $k$,
\begin{equation*}
\frac{1}{n_k} \sum_{j=1}^{n_k} h(a_j) - \int_\T h
\ge \frac{1}{4\cdot 6^k} - \frac{4}{7^k} \ge f(n_k),
\end{equation*}
which completes the proof of the theorem.


\begin{thebibliography}{100}

\bibitem{BN} P. L. Butzer and R. J. Nessel, {\it Fourier analysis
and approximation, vol. 1}, Academic Press, New York, 1971.

\bibitem{F3} A. Fujii, {\it On the uniformity of the distribution of the
zeros of the Riemann zeta function}, J. reine angew. Math. {\bf 302} 
(1978), 167--205.

\bibitem{F1} A. Fujii, {\it Some problems of Diophantine approximation
in the theory of the Riemann zeta function (III)}, 
Comment. Math. Univ. St. Pauli {\bf 42} (1993), 161--187.

\bibitem{F2} A. Fujii, {\it Some problems of Diophantine approximation
in the theory of the Riemann zeta function (IV)}, 
Comment. Math. Univ. St. Pauli {\bf 43} (1994), 217--244.

\bibitem{G} S. M. Gonek, {\it An explicit formula of Landau and
its applications to the theory of the zeta-function,}
A tribute to Emil Grosswald: number theory and related analysis,
Contemp. Math., 143, Amer. Math. Soc., Providence,
RI (1993), 395--413.

\bibitem{H} D. Hejhal, {\it On the triple correlation of the zeros
of the zeta function,} Int. Math. Res. Not. (1994), 293--302.

\bibitem{Hl} E. Hlawka, {\it \"Uber die Gleichverteilung gewisser
Folgen, welche mit den Nullstellen der Zetafunktionen zusammenh\"angen,}
Sitzungsber. \"Osterr. Akad. Wiss., Math.--Naturw. Kl. Abt. II
{\bf 184} (1975), 459--471.

\bibitem{K} J. Kaczorowski, {\it The $k$-functions in multiplicative
number theory, III. Uniform distribution of zeta zeros; discrepancy},
Acta Arith. {\bf 57} (1991), 199--210.

\bibitem{KS} N. M. Katz, P. Sarnak, {\it Zeroes of zeta functions
and symmetry,} Bull. Amer. Math. Soc. (N. S.) {\bf 36} (1999),
no 1, 1--26.

\bibitem{KR} L. Kuipers and H. Niederreiter, {\it Uniform distribution
of sequences}, New York 1974.

\bibitem{L} E. Landau. {\it \"Uber die Nullstellen der $\zeta$-Funktion},
Math. Ann. {\bf 71} (1911), 548--568.

\bibitem{M} H. L. Montgomery, {\it The pair correlation of zeros of
the zeta function,} Proc. Sym. Pure Math. {\bf 24} 
(1973), 181--193.

\bibitem{MP} M. R. Murty, A. Perelli, {\it The Pair Correlation of Zeros 
of Functions in the Selberg Class}, Int. Math. Res. Not. {\bf 10} (1999) 
531--545.

\bibitem{MZ} M. R. Murty, A. Zaharescu, {\it Explicit formulas for the 
pair correlation of zeros of functions in the Selberg class},
Forum Math. {\bf 14} (2002), no. 1, 65--83. 

\bibitem{O} A. M. Odlyzko, {\it On the distribution of
spacings between zeros of the zeta function,} Math. Comp.
{\bf 48} (1987), 273--308.

\bibitem{R} H. Rademacher, {\it Fourier analysis in number theory},
Symposium on harmonic analysis and related integral transforms: Final
 technical report, Cornell Univ., Ithica, N.Y. (1956), 25 pages;
also in: H. Rademacher, {\it Collected Works}, pp. 434--458.

\bibitem{RS} Z. Rudnick, P. Sarnak, {\it Principal
$L-$functions and Random Matrix Theory,} Duke Math.
J. {\bf 81}, {\bf 2} (1996), 269--322.

\bibitem{T} E. C. Titchmarsh, {\it The Theory of the
Riemann Zeta-function, 2nd ed. rev. by D. R. Heath-Brown,}
Clarendon Press, Oxford  1986.




\bigskip

\bigskip

\end{thebibliography}
\end{document}